\documentclass[12pt,draftcls,onecolumn,peerreview]{IEEEtran}
\usepackage[utf8]{inputenc}
\usepackage {graphicx}
\DeclareGraphicsExtensions{.eps,.pdf,.jpeg,.png}
\usepackage[cmex10]{amsmath}
\usepackage{amssymb}
\usepackage{cite}
\begin{document}
\title{Solving the Problem of Induction}
\author{Xuezhi Yang, yangxuezhi@hotmail.com}

\maketitle

\begin{abstract}
This article solves the Hume's problem of induction using a probabilistic approach. From the probabilistic perspective, the core task of induction is to estimate the probability of an event and judge the accuracy of the estimation. Following this principle, the article provides a method for calculating the confidence on a given confidence interval, and furthermore, degree of confirmation. The law of large numbers shows that as the number of experiments tends to infinity, for any small confidence interval, the confidence approaches 100\% in a probabilistic sense, thus the Hume's problem of induction is solved. The foundation of this method is the existence of probability, or in other words, the identity of physical laws. The article points out that it cannot be guaranteed that all things possess identity, but humans only concern themselves with things that possess identity, and identity is built on the foundation of pragmatism. After solving the Hum's problem, a novel demarcation of science are proposed, providing science with the legitimacy of being referred to as truth.
\end{abstract}


\section{Introduction}

Hume's problem, coined after the Scottish philosopher David Hume, represents a philosophical conundrum, which pertains to the rationality of induction and the foundation of science.

\subsection{What is Hume's Problem?}

Hume pointed out that inductive reasoning involves making predictions about future events based on past experiences. However, this kind of reasoning can not  logically  ensure its accuracy because the future may differ from the past \cite{hume1739treatise}. Hume argued that no matter how many white swans we have observed, we can not logically prove that "all swans are white" because future swans may have different colors. This raises the question that why we should believe in the effectiveness of inductive reasoning, especially considering its widespread use in scientific research.

The problem of induction is also illustrated vividly by the story of Russell's turkey. A turkey in a village observes the same thing for many days: a hand comes to feed it in the morning. For the turkey, this experience accumulates into such a strong pattern that it develops a high degree of confidence in the proposition "it will be fed every morning", until Thanksgiving arrives and shatters its belief \cite{russell1948human}.

Closely related to induction is the concept of causality. Causality implies that one event occurs as a result of another. Hume argued that people tend to believe that if one event consistently follows another, the former is the cause of the latter. However, Hume contended that this view is not drawn through logical reasoning but rather is based on habit and psychological inclinations stemming from experience.  Thus Hume believed that causality, which involves necessary connections between cause and effect, cannot be established through empirical observation or deduction.

\subsection{Responses to Hume's Problem in History}

\subsubsection{Synthetic a priori}

To address Hume's problem, Kant argued in his famous work \cite{kant1998critique} that synthetic a priori knowledge is possible. Kant reversed the empiricist programme espoused by Hume and argued that experience only comes about through the 12 a priori categories of understanding,  including the concept of causation. Causality becomes, in Kant's view, a priori and therefore universally necessary, thus sidestepping Hume's skepticism.

However, if causality cannot be derived through induction, Kant's direct assertion of causality as an a priori category may seem more brutal.

\subsubsection{Probabilistic Approaches}

The Bayes-Laplacian solution and Carnap's confirmation theory are probabilistic approaches to Hume's problem with the foundational  ideas of which are essentially the same.   

Carnap's goal is to provide a more precise measure of confirmation for scientific theories to assess their credibility \cite{1950Logical}. The measure he introduced is 

\begin{equation}
\label{Bayes}
P(h|e) = \frac{P(e|h) \cdot P(h)}{P(e)},
\end{equation}
where $e$ is the observational evidence and $h$ is a hypothesis.  A hypothesis can achieve a higher probability when supported by evidence, making it more confirmed. If the evidence is  inconsistent with a hypothesis, the degree of confirmation will be lower. Eq. (\ref{Bayes}) is actually a conditional probability expressed in the Bayesian form. 

Rule of Succession  \cite{mood1974introduction} is a formula introduced by Laplace when working on the "sunrise problem". The question here is, If event $A$ occurred  $N_A$ times in $N$ trials,  how likely would it  happen in the next trial?

Laplace considers this problem based on  such a random experiment of drawing ball randomly, with replacement, from an urn.  At each draw, a random probability $p$ is uniformly produced on the [0,1] interval, and then set $p$ portion of the balls in the urn black, and then a ball is drawn. Event $A$ is the draw is black. 

In Carnap's language, the evidence $e$= "event $A$ occurred $N_A$ times in $N$ trials", and the hypothesis $h$= "event $A$ will happen next time",  then the degree of conformation for $h$ based on $e$ is $P(h|e)$.
Laplace worked out the result as

\begin{equation}
\label{RuleofSuccession}
P(h|e)=\frac{N_A+1}{N+2},
\end{equation}
which is  widely recognized for it satisfies people's expectation. For the sunrise problem, $N_A=N$, then $P(h|e)=(N+1)/(N+2)$, which is smaller than 1 and approaches 1 when $N$ approaches infinity. This is exactly the expectation in people's mind that finite evidences don't conform a hypothesis, but more evidences improve the degree of conformation.    

Although promising and influential, Carnap's confirmation theory encountered a wide range of  criticisms from scholars among which Popper, Quine and Putnam were the most prominent figures.  The author will not elaborate on their controversies but point out the essential problem with this approach.

Carnap and Laplace's theory is based on a random experiment with an  a priori probability distribution. But are we really doing such an experiment? If so, it is easy to calculate that $e$ is an event with a very small probability. If  the sun rises with a random probability $p$ which is uniformly distributed on [0,1], it will be very rare to see 100 successive sunrises.  However, the sun has kept rising every day for millions of years.

So the essential problem with this approach is the model of the random experiment  is wrong. That's why Carnap couldn't  work out a clear and cogent theory of confirmation.
 
\subsubsection{Falsificationism}

Popper proposed "falsifiability" as the demarcation of science from pseudoscience \cite{popper1959logic}. In Popper's view, a scientific theory is a conjecture or hypothesis that can never be finally confirmed but can be falsified at any time. Propositions like "all swans are white" are not inductively derived but are conjectures. If all observed swans so far are white, the conjecture is provisionally accepted. However, once a black swan is discovered, the proposition is falsified.

From a logical perspective, if a theory $h$ implies a specific observable statement $e$, verificationism asserts that "if $e$ then $h$", which commits the fallacy of affirming the consequent. On the other hand, falsificationism claims that "if $\neg e$, then $\neg h$", which is logically sound. This logical rigor has contributed to the widespread influence of falsificationism. The idea that "science can only be falsified, not verified" has become something of a doctrine.

However, falsificationism  has its problems, with a major critique coming from Quine's holism\cite{1980From}. Holism suggests that when scientists design an experiment to test a scientific theory, they rely on various experimental equipments, which themselves embody a set of theories. For example, observing the trajectories of planets requires telescopes, which involve optical theories. In Popper's view, these theories used in experiments are treated as background knowledge and assumed to be correct. But from Quine's perspective, background knowledge can be wrong. If experimental results don't agree with theoretical predictions, it might be due to problems with the experimental equipment rather than  the theory. Thus, using experiments to falsify a theory becomes problematic.

In addition to the critique from holism, falsificationism has deeper weaknesses: how do you prove counterexamples are true? For instance, to falsify the proposition "all swans are white", you only need to find one black swan.  Assume you observed a black swan, how can you be sure that the statement "this is a black swan" is true?  Popper didn't realize  the correctness of this statement relies on the universal proposition "my senses are always correct", therefore we have no way to evade Hume and ultimately have to face the problem of induction.

This article uses modern probability theory to address the problem of induction and provides a method for calculating the confidence of a proposition over a given confidence interval, thereby solving the Hume's problem. The article emphasizes that identity is the foundation of this solution. Identity cannot be proven  but is established on the basis of pragmatism. After addressing the problem of induction, the article proposes a novel  demarcation criteria of science based on probabilistic verification.

\section{Solving the Hume's Problem}

Carnap attempted to solve the Hume's problem using a probabilistic approach, which was a step in the right direction. Unfortunately, he failed to establish a viable inductive logic  because he got the problem wrong in the first place. 

\subsection{The Concept of Probability}
From a pure mathematical perspective, probability theory is an axiomatic system. In 1933, Soviet mathematician Kolmogorov first provided the measure-theoretic definition of probability \cite{kolmogorov1933foundations}. 

Axiomatic Definition of Probability:Let $\Omega$ be the sample space of a random experiment, assign a real number $P(A)$ to every event $A$.  $P(A)$ is  the probability of event $A$ if $P(\cdot)$  satisfies the following properties:
\begin{itemize}
  \item Non-negativity: $\forall A, P(A)\geq 0$;
  \item Normalization:  $P(\Omega)=1$;
  \item Countable Additivity: If $A_m\bigcap A_n=\emptyset$ for $m\neq n$, $P\left(\bigcup\limits_{n=1}^{+\infty}A_n\right)=\sum\limits_{n=1}^{+\infty}P(A_n)$.
\end{itemize}

Although the debate on the interpretation of probability is still going on in the area of philosophy \cite{sepprobabilityinterpret}, probability theory has evolved into a rigorous branch of mathematics and has served as the foundation of information theory, which has been guiding the development of modern communication ever since its birth. 

When using the language of probability, a universal proposition is transformed into the probability of a singular proposition. For example,  "all swans are white" is expressed in probability as "the probability that a swan is white is 100\%". 

The central task of induction from the probability perspective is to estimate this probability and evaluate the accuracy of the estimation.

\subsection{How to estimate a probability?}

If the evidence $e$ is  "event $A$ occurs $N_A$ times in an  $N$-fold random experiment", our task is to estimate the probability of  event $A$, which is $p$. Maximum likelihood algorithm has been a well accepted criterion for optimal estimation. In this criterion, we choose the estimation to maximize the probability of  $e$,  which is expressed as 

\begin{equation}
P(e)=\begin{pmatrix}
       N   \\
       N_A
\end{pmatrix}p^{N_A}(1-p)^{N-N_A}.
\end{equation}
Let
\begin{equation}
\frac{dP(e)}{dp}=0,
\end{equation}
we get the maximum likelihood estimation of $p$, denoted as
\begin{equation}
\hat{p}=\frac{N_A}{N}.
\end{equation}

To assess the accuracy of this estimation, we need the concept of confidence on confidence interval.

\subsection{Confidence on Confidence Interval}

Because $p$ is a real number, there is no chance the estimation is exactly the same as $p$. Therefore, confidence should be defined on a confidence interval with a non-zero width. Confidence is the probability that the true value falls within this interval. 

From the perspective of maximum likelihood estimation, the optimization goal of this problem should be, given a width of the confidence interval, find a confidence interval  that maximize the confidence. Earlier, we have obtained the maximum likelihood estimation of $p$ as $N_A/N$. In this paper, we just simplify this problem to set the confidence interval as $D$, so that $N_A/N \in D$. 

To calculate the confidence, we first need to make an assumption that in absence of any observed facts, $p$ is uniformly distributed on $[0, 1]$.  Under this assumption, for each possible probability, we calculate the probability of $e$ to form a curve. The ratio of the area under the curve within $D$ to the total area under the curve is the probability that the true value falls within $D$, which is the confidence. Then, the confidence is given by

\begin{equation}
c=\frac{\int_{D}x^{N_A}(1-x)^{N-N_A}dx}{\int_0^1 x^{N_A}(1-x)^{N-N_A}dx}.
\end{equation}

Although the uniform prior over $p$ is the same as the Principle of Indifference in Rule of Succession,  the ideology is different. Unlike Laplace's solution, the random experiment in this approach is, event $A$ happens with a fixed probability rather than a random one.  It is just that, we don't know the true value of $p$. The task here is to use evidence $e$ to estimate $p$. The uniform prior is used here to define confidence, which is  a subjective criterion of evaluation. In contrast, it is constitutive of the random experiment in Laplace's approach \cite{sepinductionproblem}. 

Let's take a simple example. Based on the fact that "all $N$ swans were observed to be white", calculate the confidence  of "the probability of a white swan being greater than 90\%", where the confidence interval is $[0.9, 1]$. The expression for this confidence  is

\begin{equation}
c=\frac{\int_{0.9}^1x^Ndx}{\int_0^1 x^Ndx}=1-0.9^{N+1}.
\end{equation}

\begin{figure}[h]
\begin{center}
\includegraphics[width=3in]{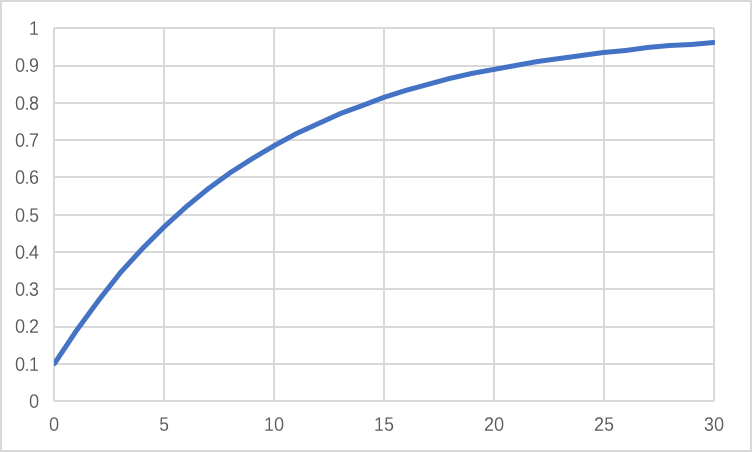}
\caption{Confidence on $[0.9, 1]$ with respect to $N$}\label{whihtegoose}
\end{center}
\end{figure}

As shown in Fig. 1,  it can be seen that when $N=10$, the confidence is approximately 68.6\%. When $N=30$, the confidence increases to 96.2\%. If the confidence interval is reduced to $[0.99, 1]$, then $N=300$ is required to achieve a confidence  of 95.1\%.

With the concepts of confidence  and confidence interval, we can make a response to Russell's turkey. After two months of feeding by the farmer, the conclusion drawn by this turkey should have been "the confidence  of 'the probability of feeding is greater than 99\%' is $1- 0.99^{61}=45.8\%$", so it is not surprising that it was killed by the farmer on Thanksgiving.  Besides, for a cow who has lived on the farm for 10 years, it is normal to conclude that a turkey has a high possibility to be slaughtered on Thanksgiving.

A person around 30 years old who has seen the sun rise 10000 times can conclude that the confidence of "the probability of  sunrise is greater than 99.9\%" is $1-0.999^{10000}=99.995\%$. If he also believes in historical records that humans have seen sunrise for over a million years, then the confidence interval can continue to narrow and the confidence is closer to 100\%.

\subsection{The Law of Large Numbers}

The law of large numbers is the fundamental law of probability theory, which comes in various forms, such as Bernoulli, Sinchin, Chebyshev's law of large numbers, and the strong law of large numbers. Let's take a look at the most fundamental form, Bernoulli's law of large numbers.

Bernoulli's Law of Large Numbers: If a random experiment is conducted $N$ times, event $A$ occurs $N_A$ times, and $p$ is the probability of $A$, then for any $\varepsilon > 0$

\begin{equation}
\displaystyle \lim_{N\rightarrow \infty}P\left(\left|\frac{N_A}{N}-p\right|< \varepsilon\right)=1.
\end{equation}

The proof of this theorem can be found in any textbook of probability theory. Bernoulli's theorem of large numbers states that when the number of trials is sufficiently large, the value of $N_A/N$ approaches $p$ infinitely in probability.

From another perspective, $P\left(\left|N_A/N-p\right|< \varepsilon\right)$ is the confidence of $p$ being located in the confidence interval $[N_A/N-\varepsilon,N_A/N+\varepsilon]$. Bernoulli's law of large numbers states that for any small confidence interval, when $N $ is large enough, the confidence can be infinitely close to 1.

\subsection{Degree of Conformation}

We have discussed the concept of  confidence on confidence interval. If the two parameters still seem too complicated, we can further simplify them to one parameter. Suppose the width of the confidence interval is $d$, then we can introduce a metric, degree of conformation,  as 

\begin{equation}
C=\max_d (1-d)\cdot c.
\end{equation}

For the sunrise problem, $N_A=N$, and  the confidence on $[x,1]$ is $1-x^{N+1}$, By simple inference we get 

\begin{equation}
C=\frac{1}{\sqrt[N+1]{N+2}}\frac{N+1}{N+2}.
\end{equation}

When $N=2$, $C=0.47247$, and when $N=10000$, $C=0.9990$, while the results of Rule of Succession are 0.75 and  0.9999, respectively . Though it is of little meaning to compare the detailed numbers, the conformation provided by this solution is more conservative than that of the Rule of Succession.

\section{Identity}

The basis for solving the Hume's problem is the existence of probability. So let's further ask, does probability exist? This involves the concept of identity. 

The so-called identity means invariance. Only things with identity can humans understand and use the acquired knowledge about them to guide practices. If a thing changes after people understand it, then the previous understanding becomes useless.

We acknowledge that there are things in the real world that do not have identity. For example, a shooting star in the night sky disappears after a few seconds of brilliance. You see this scene and tell your friend, 'Look, meteor!'. When your friend looks up, it's already gone. Although your friend may not think you are lying, he has no evidence to  acknowledge you because this phenomenon is non-repeatable, and no one will care anymore after its disappearence.

What humans are concerned about are things with identity. For example, if you see an apple in front of you,  you look at it once, you look at it again, you look at it ten times, it will keep being an apple in your eyes. If you let your friend look at it, he will also see an apple. So, the existence of the apple has identity. For this apple, you can say to your friend, "Do you want to eat this apple?", your friend might say, "That's great!" and then pick it up and eat it. Humans are able to communicate and cooperate on objects with identity.

Hume's questioning of induction and causality is essentially a negation of identity. Russell's turkey was killed on Thanksgiving, so how can we guarantee that the apple in front of us won't suddenly turn into a rabbit? We successfully responded to this challenge using a probabilistic approach. That is to say, our assumption of identity is not that the apple has been and will always be an apple, but rather that there is a probability that it is an apple,  without limiting this probability to be 100\%, which preserves the possibility of "an apple suddenly becomes a rabbit" to avoid dogmatism to the utmost extent. 

The notion of identity assumes there is an invariant probability, which is the basis for our argument for all other propositions. It is a presupposition and cannot be proven because there is no more fundamental propositions to prove it.

\section{Practicality is the Foundation of Identity}

If  $p$ exists, then according to the law of large numbers, $N_ A/N$ approaches $p$ infinitely after $N$ is sufficiently large. If we conduct experiments, such as flipping a coin, we do observe a phenomenon where the ratio of the number of heads ($N_A $) to the total number  ($N $) gradually approaches a constant value. However, this does not prove the existence of probability. Because according to Hume's query, although  $N_ A/N$ tends to be a constant value,  it is only a case of finite number of trials and cannot be extrapolated to an infinite conclusion. This is also the reason why we say identity cannot be proven.

Let's do a thought experiment. Suppose the demon of Descartes want to subvert our belief in identity and manipulate our experiment of flipping coins. The devil's purpose is to make $N_ A/N$ has no limit by letting $N_ A/N$ oscillates in constant amplitude.

The devil first makes the coin flip normally, so $N_ A/N $ will be around 0.5, and then manipulates the probability of  heads be 0.6 until $N_ A/N$ reaches 0.55, and then  manipulates the probability of  heads be 0.4 until $N_A/N$ reaches 0.45 and cycles like this, then $N_ A/N$ oscillates between 0.45 and 0.55 without convergence, meaning there is no limit for  $N_ A/N$. So can the demon's approach overturn our belief in the existence of probability?

However, it cannot. Because $N$ increases with each cycle, so changing the overall statistical characteristics in the next cycle will require more trials, which results in an exponential increase in the number of trials per cycle. At the beginning, due to the short cycle, the results are not reliable, so these experimental results could not guide human activities and no one would care. Later on, as the cycle becomes longer, for those who make short-term predictions, they will find that the current probability is 0.6 or 0.4, and they can arrange their practical activities based on this result and achieve success. It cannot be ruled out that all the physical laws currently mastered by humans are only invariants arranged by demons within a cycle. For long-term observers, it is easy to detect the oscillation pattern of  $N_ A/N$, which is a broader sense of identity.

This is actually the same as how we handle our daily lives. We believe an apple remains unchanged in the short term, and this identity can guide short-term activities, such as dealing with questions such as "Do you want to eat the apple?". But in the long cycle, apples will go through the process of  freshness, loss of luster, wrinkling, decay and blackening. Although the apple has changed, the process of change follows the same pattern and remains unchanged, which is a long-term identity. Based on this identity, humans can derive the most favorable principles of action for themselves, such as eating apples before they become stale.

Of course, the demon can also come up with more complex ways to subvert the assumption of identity, but this is not of much use. Because in the face of a constantly changing world, humans always extract those unchanging characteristics, identify their patterns, and use them to guide their practical activities. For those changes that cannot be mastered, the human approach is to ignore them and let them go. 

So, the assumption of identity is based on pragmatism, which is the final foundation of human knowledge \cite {Yang2023}.  

\section{Demarcation of Science}

The demarcation problem is the fundamental issue in the philosophy of science. The most influential demarcation criterion is Popper's falsification principle. According to Popper's theory, a proposition is  scientific  if it has falsifiability.

The problem with Popper's  standard is that it is too loose so that many nonsense propositions will also be considered scientific. Such as astrology, which was rightly taken by Popper as an example of pseudoscience, DOES have falsifiability, and has in fact been tested and refuted.

Popper's falsification principle is a stopgap measure before the Hume's problem is solved,  after which we can work out a precise definition of science.

\begin{figure}[h]
\begin{center}
\includegraphics[width=3in]{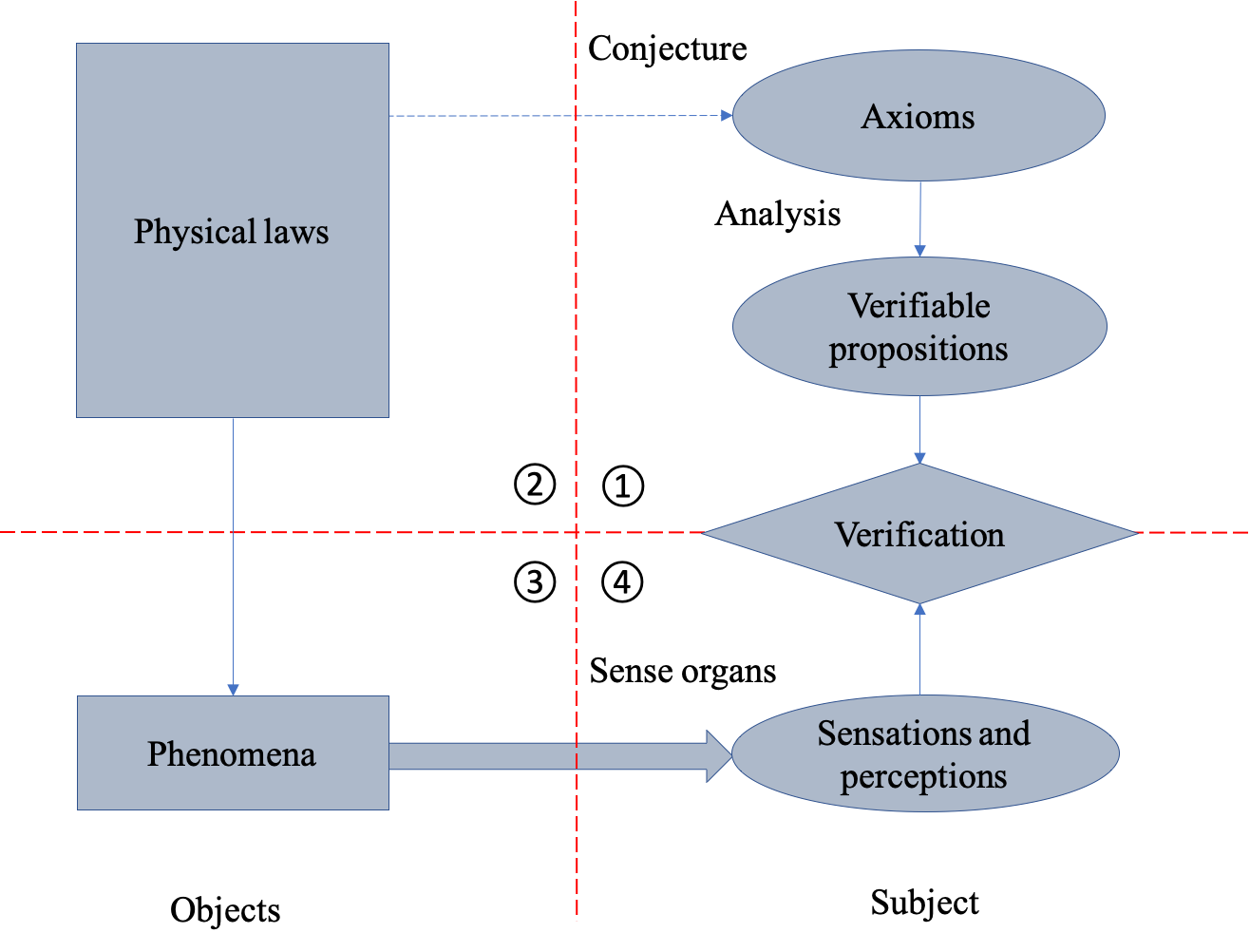}
\caption{Definition of science}\label{science}
\end{center}
\end{figure}

The definition of science is shown in Figure \ref{science}. The left half plane represents the part of the world with identity, while the right half plane represents the subject, which is human. 

The third quadrant is miscellaneous phenomena that trigger our sensations and perceptions, while behind the phenomenon are  physical laws, located in the second quadrant. Plato believed that there exists a world of ideas, and objects in the real world are copies of ideas. In our model, there exists a world of physical laws. Unlike Plato, physical laws are not blueprints of phenomena, but rather to constrain and regulate them.

Popper believes scientific theories are conjectures of scientists on what laws are, which is in the first quadrant.  Induction is not a logical method, but a way of conjecture. The results of conjecture appear in the form of axioms, becoming the starting point of reasoning.  Conjectures should obey the Occam's Razor Principle, also known as the economic principle of thinking, by cutting off the unverifiable parts of them.

Axioms are usually universal propositions and cannot be directly verified. Then it is necessary to combine actual scenarios to deduce a proposition that can be empirically confirmed or falsified. For example, we can conjecture  "the sun rises every day" and use this proposition as an axiom. This axiom cannot be directly verified, and verifiable propositions need to be derived from it through logic, such as "the sun rose yesterday", "the sun will rise tomorrow", "the sun will rise the day after tomorrow", and so on. If all and large enough amount of propositions derived from this axiom are verified to be correct, then "the sun rises every day"  is confirmed with high degree of confirmation.

Newton's laws, theory of relativity, and quantum mechanics are also judged according to this standard. For example, Newton's laws have many application scenarios, such as free falling bodies, the trajectory of cannonballs and planets, and so on. In each scenario, a series of verifiable propositions are derived. If these propositions are experimentally confirmed, then Newton's laws are verified with high conformation  in these scenarios. 

In scenarios in which Newton's laws are yet to be verified, we can conjecture that they will be still valid. But experiments have shown that Newton's laws do not hold true in scenarios close to the speed of light and atomic scale, so are falsified in these scenarios. But that does not affect conclusions in low-speed macroscopic scenarios. Once a theory is conformed with high degree of confirmation in a scenario, it is only theoretically possible and practically impossible to falsify it again under that scenario. Therefore, Newton's laws deserve the title of TRUTH in confirmed scenarios. 

From the above discussions, we can summarize the four elements of science: Conjecture, Logic, empirical Verification, and Economics. The use of probability methods to address the Hume's problem is reflected in the verification part.

This standard of demarcation is actually the revival and development of logical positivism or logical empiricism, and is much stricter than the falsification criterion. Probabilistic verification requires a large number of experiments to obtain a high degree of confirmation, so any establishment on luckiness cannot pass the verification criteria.  In addition, according to this standard, not only Newton's laws, theory of relativity, and quantum mechanics are science, but common senses of life such as "the sun rises every day" and "apple is eatable" are also science. This way, the notion of science will  enter the everyday lives of ordinary people  and can play a role in improving the scientific literacy of the whole people.

\section{Conclusions}
 
The Hume's problem is a fundamental problem that runs through epistemology. Kant and Popper were unable to solve this problem and adopted an evasive attitude. Carnap's attempt was an endeavor in the right direction, but ended up in failure. The goal achieved in this article is exactly what Carnap  attempted to achieve. For three hundred years, the Hume's problem has been a dark cloud over science and relegated science to the position of "can never be confirmed and wait to be falsified". After the solution of the Hume problem, science is confirmed to deserve the title of truth in a probabilistic sense.

\bibliographystyle{IEEEtran}
\bibliography{IEEEfull,philosophy}

\end{document}